\documentclass[12pt]{amsart}
\usepackage{a4}
\usepackage{amssymb}
\usepackage{amsmath}
\usepackage{amsthm}
\usepackage{amstext}
\usepackage{amscd}
\usepackage{latexsym}
\usepackage{graphics}
\usepackage{color}

\swapnumbers
\theoremstyle{plain}
\newtheorem{thm}{Theorem}[section]
\newtheorem{q}[thm]{Question}
\newtheorem{lemma}[thm]{Lemma}

\newtheorem{cor}[thm]{Corollary}
\newtheorem{corollary}[thm]{Corollary}
\newtheorem{prop}[thm]{Proposition}

\theoremstyle{definition}
\newtheorem{rem}[thm]{Remark}

\newtheorem{rmk}[thm]{Remark}

\newtheorem{defn}[thm]{Definition}
\newtheorem{ex}[thm]{Example}

\newtheorem{notation}[thm]{Notation}

\newcommand{\rdown}[1]{\lfloor{#1}\rfloor}

\newcommand{\supp}{{\rm Supp}}

\newcommand{\rank}{{\rm rank}}
\newcommand{\Ker}{{\rm Ker \ }}
\newcommand{\Pic}{{\rm Pic}}

\newcommand{\Aut}{{\rm Aut}}

\newcommand{\Tr}{{\rm Tr}}

\newcommand{\PGL}{{\rm PGL}\,}   
      
\newcommand{\GL}{{\rm GL}}
\newcommand{\SL}{{\rm SL}}
\newcommand{\PSL}{{\rm PSL}}

\newcommand{\Hes}{{\it Hes} }

\renewcommand{\tilde}{\widetilde}

\newcommand{\sO}{{\mathcal O}}

\newcommand{\C}{{\mathbb C}}

\newcommand{\F}{{\mathbb F}}

\newcommand{\N}{{\mathbb N}}

\renewcommand{\P}{{\mathbb P}}
\newcommand{\Q}{{\mathbb Q}}

\newcommand{\Z}{{\mathbb Z}}

\input{xy}
\xyoption{all}

\begin{document}
\title{Notes on $\pi_1$ of Smooth Loci of Log Del Pezzo Surfaces}
\author{Chenyang Xu}

\email{chenyang@math.princeton.edu}
\address{Department of Mathematics, Princeton University, Princeton NJ 08544.}

\date{Oct. 26, 2008}

\begin{abstract}
It is known that the fundamental groups of smooth loci of log del Pezzo Surfaces are finite groups. The aim of this note is to study these finite groups. A short table containing these groups is given. And lots of groups on the table are proved to be fundamental groups.
\end{abstract}

\maketitle
\vspace*{6pt}
\section{Introduction}
A projective surface $R$ over $\C$ is called a {\it log del Pezzo surface}, if it contains only quotient singularities, and the canonical divisor $K_R$ is an anti-ample $\Q$-divisor. Although the fundamental group of $R$ is always trivial, the fundamental group of the smooth locus $\pi_1(R^{sm})$ is in general not zero. Nevertheless, it is known such a group is always finite (cf. \cite{zhang}, \cite{keel}). The aim of this paper is to determine these groups. 

Our approach to this problem is as follows. Given a log del Pezzo surface $R$, we take the universal cover of its smooth locus $R^{sm}$.  Knowing $\pi_1(R^{sm})$ is finite (cf. \cite{zhang}, \cite{keel}), the Riemann Existence Theorem (cf. \cite{sga1}) says that the universal cover is actually an algebraic variety. Therefore, we can take the normal closure $S$ of $R$ in the function field of this covering space. In this way, we get a pair $(S,\pi_1(R^{sm}))$, where $S$ is also a log del Pezzo surface, and $\pi_1(R^{sm})$ is a finite group acting on it, such that for every nontrivial element $g\in \pi_1(R^{sm})$, the fixed locus $S^g$ is isolated. We can also equivariantly resolve $S$ to get a smooth rational surface carrying the same finite group action. This motivates the following definitions,

\begin{defn} 
We call a finite group $G$ acting on a normal surface $S$ {\it an action with isolated fixed points} (IFP), if $S$ has at worst quotient singularities, and for every nonunit element $g\in G$, the fixed locus $S^g$ consists of finite points. Similarly, we call $(S,G)$ {\it birational to an action with} IFP if there is a $G$-equivariant birational proper model $S'$ of $S$, such that $(S',G)$ is an action with IFP. 
\end{defn}

Now we can divide our question into 3 parts: 
\begin{enumerate}
\item[(1)] finding all the birational classes $(S,G)$ containing a representative $(\tilde{S},G)$ with IFP;
\item[(2)] determining those groups $G$, for which we can choose  $(\tilde{S},G)$ as in (1) with the additional property that $K_{\tilde{S}}$ is anti-ample; and
\item[(3)] for any $G$ appearing in the Step (2), checking the existence of $(\tilde{S},G)$ satisfying $\pi_1(\tilde{S}^{sm})=e$.
\end{enumerate} 

In a recent paper \cite{dol}, all finite subgroups of the Cremona group are classified. Based on their table, we can solve the Problem (1).
\begin{thm}\label{main}
Let $G$ be a finite group which acts on a rational surface $S$ such that $(S,G)$ is birational to an action with {\rm IFP}, then $G$ precisely is one of the following groups:
\begin{enumerate}
\item a finite subgroup $G$ of $\GL_2(\C)$ whose abelian subgroups are all cyclic,
\item a finite subgroup $G$ of $\PGL_2(\C)\times \PGL_2(\C)$ whose subgroups as $G_1 \times G_2$ have the property that $|G_1|$ and $|G_2|$ are coprime,
\item $\Z/n:\Z/3$ or $\Z/2\times (\Z/n:\Z/3)$, where $n$ is an odd integer and $\Z/n:\Z/3$ means the group generated by $u:(x_0,x_1,x_2)\to (x_1,x_2,x_0)$ and $v:(x_0,x_1,x_2)\to (\epsilon_n x_0,\epsilon_n^s x_1,x_2)$ $(s^2-s+1\equiv 0 \mod n)$ where $\epsilon_n$ is a primitive n-th root,
\item $F_{4n}$, $G_{4n}$ and $H_{4n}$ which are groups of order $4n$ (for the definition, see (\ref{f_0})), or 
\item $(\Z/3)^2:\Z/2$, $(\Z/3)^2:\Z/4$ and $(\Z/3)^2:Q_8$.
\end{enumerate}
\end{thm}

Although Theorem(\ref{main}) is the strongest statement, we emphasize that there is a more conceptual version as follows:
\begin{thm}
Given a finite group $G$ which acts on a smooth projective rational surface $S$, it is an action birationally with IFP if and only if it satisfies the following conditions:
\begin{enumerate}
\item[(1)] for any point $x\in S$, every abelian subgroup of the stabilier $G_x$ is cyclic, and
\item[(2)] for any nonunit element $g\in G$, every curve $C\subset S^g$ satisfies genus $g(C)=0$.
\end{enumerate}
\end{thm}

In fact, in our case by case study, the phenomena can be formulated as a dictonomy: 
\begin{enumerate}
\item for a simple rational surface, i.e. a Hirzebruch surface or a del Pezzo surface of degree $\ge 5$, the minimal action is birationally with IFP, if and only if for any point $x$, the stabilizer $G_x$ does not contain a noncyclic abelian subgroup;
 \item for other complicated rational surface, any minimal action contains a nontrivial element, which fixes a positive genus curve pointwise. 
\end{enumerate}

In \cite{kollar}, a similar method is used to study the case when $G$ is abelian, yielding a list of possible first homology groups of log del Pezzo surfaces. By listing the abelian groups in the above table, we can refine the results there (cf. \cite{kollar}, $11$).  
\begin{corollary}
Let $S$ be a log del Pezzo surface. Then
$H_1(S^{sm}, \Z)$ is one of the following groups: $(\Z/3)^2$, $Z/3\times \Z/6$, $\Z/2\times \Z/n$ ($n$ is $4$ or $4k+2$) or $\Z/m$ for any $m$.
\end{corollary}

Then applying the equivariant minimal model program, we can also answer the question (2), namely
\begin{thm}\label{main2}
If we can choose $(\tilde{S},G)$ in Theorem(\ref{main}) satisfying the additional property: $\tilde{S}$ is a log del Pezzo surface. Then $G$ is precisely one of the groups listed in (1)-(4) there.
\end{thm}

In the last section, we aim to solve the third problem. We construct models $(\tilde{S},G)$ which satisfy the property $\pi_1(\tilde{S}^{sm})=e$ for most groups $G$ in (\ref{main}). Unfortunately, we leave three series of groups undetermined.

\noindent{\bf Acknowledgement}: I am indebted to my advisor, J\'anos Koll\'ar, for suggesting me this question, also for his useful comments and crucial suggestions. Without his encouragement, this paper would have never been written up. I thank Igor Dolgachev for sending me the new version of [DI06]. As mentioned above, our work relies on it substantially. I also thank Ruochuan Liu, Garving Luli, Yi Ni and Zhiwei Yun for helpful conversations. 

\section{Local and Global Results}

For a table of the quotient surface singularities, see (\cite{Br}). Notice that any quotient singularity is rational. In particular, given a resolution, the irreducible components of the exceptional locus are smooth rational curves, and the configuration of the exceptional locus is a tree.

Now let $S$ be a projective rational surface with $G \subset \Aut(X)$ a finite group of automorphism. In this section, we aim to give some criteria to decide whether $(S,G)$ is birational to an action with IFP. By abuse of notation, we use $\sum_{S,G}$ to denote both the set of the irreducible curves which are fixed pointwise by some nontrivial elements of $G$ and the corresponding reduced divisor. When the surface and the group are clear, we will omit the subscript. We call a subset $\{C_1,C_2,\cdots,C_n\}\subset\sum_{S,G}$ a $cycle$ if the intersection of $C_i$ and $C_{i+1}$ are nonempty for all $1\le i \le n$ ($C_{n+1}=C_1$). We also use the same letter to mean both the divisor and its birational transforms on any other birational models. We also define a subset 
$$\tilde{\sum}=\{C\in \sum|C\mbox{ intersects with other curves of } \sum \mbox{at most at two points.} \}$$

Let $x$ be a point in the smooth $G$-surface $S$, with the stabilizer $G_x \subset G$. In an analytical neighborhood of $x$, it is isomorphic to a germ $(\C^2,0)$ with a finite group $G_x\subset \GL_2$ action. For the lemma below, we assume this identification.

\begin{lemma}\label{noncyclic}
If $G_x$ contains a noncyclic abelian group, then there are two curves $C_1, C_2\in \sum$, such that if we denote $\langle g_1 \rangle$ (resp. $\langle g_2 \rangle$) the subgroup fixing $C_1$ (resp. $C_2$) pointwise, then $\langle g_1, g_2 \rangle \subset G_x$ gives a noncyclic abelian group. Furthermore, for any sequence of $G_x$-equivariant blow-ups $\pi: S \to \C^2$, $C_1, C_2$ belong to the same connected component of $\sum_S$.
\end{lemma}
\begin{proof} Since $G_x$ contains a noncyclic abelian group, we know $(\Z/p)^2\subset G_x$ for some prime $p$. Given a $G$-surface $S$, $\sum_{(\Z/p)^2} \subset \sum_{G}$. So it suffices to prove for $G_x\cong (\Z/p)^2$. But if we choose the basis suitably,  any embedding of $(\Z/p)^2$ to $\GL_2(\C)$ is conjugate to the embedding  $(0,1) \to diag\{1,e^{\frac{2\pi i}{p}}\}$ and $(1,0) \to diag\{e^{\frac{2\pi i }{p}},1\}$, hence the first statement is obvious.

To prove the second statement, we can just look at the blow-ups whose centers are the intersection points of (at least) two curves in $\sum_S$. Assume after blowing up $x$, we get $\pi : S_1 \to \C^2$, with the exceptional divisor $E=\pi^{-1}(0)$.
The action can be lifted from $\C^2$ to $S_1$ such that $G_x$ acts on $E$ through the natural homomorphism $p:\GL_2(\C) \to \PGL_2(\C)$. Its restriction on $(\Z/p)^2$ has a nontrivial kernel which fixes $E$ pointwise. 

Furthermore, $(\Z/p)^2$ fixes the two intersection points of $C_i (i=1,2)$ and $E$, so if we replace $C_1$ and $C_2$ by $C_i$ and $E$, the assumptions of the lemma still hold. Repeating the argument, we can see there does not exist any sequence of equivariant blow-ups $\pi: S \to \C^2$ which separates the birational transforms of $C_1$ and $C_2$ into different connected components of $\sum_S$.  
\end{proof}

This local computation leads to a global version:
\begin{cor}\label{criterion}
Let a finite group $G$ act on a smooth surface $S$. Assume $\sum_S$ contains a cycle $\{C_1,C_2,\cdots,C_n\}$ and $\langle g_i \rangle$ fixes $C_i$ pointwise. If for every $i$, $\langle g_i, g_{i+1} \rangle \subset G$ is noncyclic abelian, then the action is not birational to one with only {\rm IFP}. 
\end{cor}

We can also prove the converse of (\ref{noncyclic}), namely
\begin{lemma}\label{separation}
Notations as (\ref{noncyclic}), if $G_x$ is cyclic, and $C_1,C_2$ are the curves fixed by some nontivial subgroups of $G_x$ (there are at most 2 such curves), then there is a sequence of $G_x$-equivariant blow-ups $\pi:S\to \C^2$, such that the strict transforms of $C_1,C_2$ are disconnected in the configuration of $\sum_S$. Given a $C_i$, we can assume that all the exceptional curves $E_i\subset \sum_S$ is a tail added to it. 
\end{lemma}
\begin{proof} We can assume the order-$r$ cyclic group $G_x$ is generated by
$(e^{\frac{2p\pi i}{r}}, e^{\frac{2q\pi i}{r}})$, where $(p,q)=1$.  Blowing up $x$, we will change an intersection point from type $\frac{1}{r}(p,q)$ to two intersection points of types $\frac{1}{r}(p,q-p)$
and $\frac{1}{r}(p-q,q)$. Keep blowing up the new intersection points and changing the action in this way, we can choose $t$ such that $(r,(t+1)p-q)=1$, so the image of $\frac{1}{r}(p,q-tp)$ in $\PGL_2(\C)$ is of order $r$. This means the only element of $G_x$ which fixes the corresponding exceptional curve is the unit. 
\end{proof}

The global version of (\ref{separation}) is, 
\begin{lemma}\label{cycstab}
If a finite group $G$ acts on a smooth rational surface $S$ with the following properties:
\begin{enumerate}
\item for any nonunit element, its fixed locus only consists of smooth rational curves and isolated points;
\item every intersection point of 2 curves in $\sum_S$ has an abelian stabilizer; and
\item after separating the intersection points contained in $\tilde{\sum}$ and with cyclic stabilizer, every component of the configuration of $\sum_S$ is a chain, 
\end{enumerate} 
then $(S,G)$ is birational to an action with {\rm IFP}.
\end{lemma}

\begin{proof}
For every point $x$ satisfying condition (3), by the previous lemma, we know there exists a sequence of $G$-equivariant blow-ups $S'\to S$ which separates the two branches in $\sum_S$ containing $x$ and only adds a tail to the component in $\tilde{\sum}$. Hence, the configuration of $\sum_{S'}$ is a disjoint union of chains.

Then to contract $\sum_{S'}$, we want that the self-intersection of each component in $\sum_{S'}$ is less or equal to $-2$. This may not be true for $S'$. However, we can blow-up general orbits on curves in $\sum_{S'}$. The exceptional locus we create in this way satisfies the property that for any nontrivial element $g\in G$ acting on it, the fixed locus is isolated. Hence after a sequence of such blow-ups, we can assume there is a surface $S''$ with a proper $G-$equivariant birational morphism $f':S''\to S'$ inducing an isomorphism from $\sum_{S''}$ to $\sum_{S'}$, and all the irreducible components in $\sum_{S''}$ have self-intersection numbers less or equal to -2. From the table of quotient surface singularities  (cf. \cite{Br}), we can contract every connected components of $\sum_{S''}$ (which is a chain) to a cyclic quotient surface singularity. 
\end{proof}

\begin{lemma}\label{abelstab}
For any $(S,G)$, we have a equivariant blow-up $\pi:S'\to S$, such that the condtion(2) of the above lemma holds for $S'$. 
\end{lemma}
\begin{proof}
Assume we have a point $x\in S$ with a nonabelian stabilizer $G_x\in \GL_2(C)$. Blow up $x$, we have an exceptional divisor $E$. Then for any point $y\in E$, the stabilizer $G_y$ fits the exact sequence
$$1\to K \to G_y \to H\to 1,$$
where $K$ is the kernel of $\rho:G_x \to \PGL_2(\C)$, and $H$ is the stabilizer of $y$ for the induced action $im(\rho)$ on $E$. Since $H$ is abelian, and $K$ is in the center of $\GL_2(\C)$, we conclude that $G_y$ is an abelian group.  
\end{proof}

\noindent A priori, (\ref{criterion}) and (\ref{cycstab}) do not summarize all possible cases, but together with the following simple lemma, it is enough for our purpose.

\begin{lemma}\label{highgenus}
Let $g$ be a nontrivial element in $G$, if there is a curve $C\in S^g$ with genus $g(C)> 0$, then $(S,G)$ is not birational to an action with {\rm IFP}.
\end{lemma}
\begin{proof}
For any birational $G$-map $f:S\dashrightarrow S'$, if $S'$ has only quotient singularities, then $f$ cannot contract $C$, so $(S',G)$ cannot be with only IFP.  
\end{proof}

\section{Birational Classification of Actions on Rational Surfaces with Isolated Fixed Points }

The aim of this section is to give a complete classification of actions on a rational surface $(S,G)$, which have birational $G$-models with only IFP. For a given surface $S$, the actions are classified up to the conjugation of the {\it automorphism} of $S$. Since it is a birational property, we will only consider the minimal actions, i.e, if there is a birational $G$-morphism  $f: S \to S'$, then it is in fact an isomorphism.

\begin{notation}\label{notation}
We employ some standard notations for groups here:
\begin{enumerate}
\item[$\bullet$] $\Z/n$ means the order-$n$  cyclic group;
\item[$\bullet$] $S_n$, the {\it permutation group} of degree $n$;
\item[$\bullet$] $A_n$, the {\it alternating group} of degree $n$;
\item[$\bullet$] $D_{2n}$, the {\it dihedral group} of order $2n$;
\item[$\bullet$] $Q_{4n}=\langle  a, b | a_{2n} = 1, b_2 = a_n, b^{-1}ab = a^{-1}\rangle$, {\it dicyclic group} of order $4n$, a {\it generalized quaternion group} if n = 2k; 
\item[$\bullet$] $L_n(q) = \PSL(n, \F_q)$, where $q = p^r$ is a power of a prime number $p$; 
\item [$\bullet$] $H_n(p)$, the Heisenberg group of unipotent $n \times n$-matrices with entries in $\F_p$;
\item[$\bullet$] $A\bullet B$ is an upward extension of $B$ with help of a normal subgroup $A$;
\item[$\bullet$] $A : B$ is a split extension, i.e. a semi-direct product $A \rtimes B$ (it is defined
by a homomorphism $\varphi : B \to \Aut(A)$);
\item[$\bullet$]  $A\wr S_n$ is the wreath product, i.e. $A^n:S_n$ and $S_n$  acts
on $A^n$ by permuting the factors;
\item[$\bullet$] $(G_1,H_1,G_2,H_2)_{\alpha}$ means the subgroup of $G_1 \times G_2$ consisting of elements
$\{(g_1, g_2)|g_1$ and $g_2$ has the same image under the isomorphism $\alpha : G_1/H_1 \cong G_2/H_2\}$. We will omit $\alpha$ if the isomorphism is clear;
\item[$\bullet$] $\mu_n$, the group of $nth$ roots of unity with generator $\epsilon_n = e^{2\pi i/n}$.
\end{enumerate}

We also need the notations for polyhedron groups which are precisely all possible finite subgroups of $\PGL_2$:
\begin{enumerate}
 \item [$\bullet$] a cyclic group $\Z/n$ of order $n$;
\item [$\bullet$] a dihedral group $D_2n$ of order $2n$;
\item [$\bullet$] the tetrahedral group $T \cong A_4$ of order 12;
\item [$\bullet$] the octahedral group $O \cong S_4$ of order 24;
\item [$\bullet$] the icosahedral group $I \cong A_5$ of order 60.
\end{enumerate}
We will use $\bar{T}, \bar{O}$ and $\bar{I}$ to mean the corresponding double cover of $T,O$ and $I$ under the homomorpshim 
$$\SL_2\to \PGL_2.$$
Notice that the double cover of $D_{2n}$ is $Q_{4n}$.
\end{notation}

\begin{rmk}
For the following discussions in this section, we will heavily rely on the results in \cite{dol}. In fact, we will do a case by case study for Section 4-6 of their paper.    
\end{rmk}

\noindent {\bf The Case} $\mathbf{S=\P^2}:$

Recall some standard terminology from the theory of linear groups. Let $G$ be a finite subgroup of the general linear group $\GL(V)$ of a complex vector space $V$. The
group $G$ is called {\it intransitive} if the representation of $G$ in $V$ is reducible. Otherwise
it is called {\it transitive}. A transitive group $G$ is called {\it imprimitive} if it contains an
intransitive normal subgroup $G'$. In this case $V$ decomposes into a direct sum of
$G'$-invariant proper subspaces, and elements from $G$ permute them. A group is
{\it primitive} if it is neither intransitive, nor imprimitive. We reserve this terminology
for subgroups of $\PGL(V)$ keeping in mind that each such group can be represented
by a subgroup of $\GL(V)$. 

As any element $g$ of finite order in $\PGL_3(\C)$ can be lifted as an element of $\GL_3(\C)$, $g$ fixes a curve pointwise if and only if the characteristic polynomial of the lifting has multiple roots, in which case, the curve is a line.\\ 

{\it Intransitive actions:} for any intransitve group action $(\P^2,G)$, $G$ also linearly acts on $\C^2$ with an equivariant embedding $i: \C^2 \to \P^2$, so $G\subset \GL_2$. 

\begin{prop}\label{intransitive}
An intransitive action $G$ on $\P^2$ is birational to an action with {\rm IFP} if and only if any abelian subgroup $H \subset G$ is cyclic.  
\end{prop}
\begin{proof} 
Blowing up the origin of $\C^2$, we know if any abelian subgroup of $G$ is cyclic, then the conditions of (\ref{cycstab}) all hold for this ruled surface.  To prove the ``only if" part, we notice that if $G$ has a noncyclic abelian subgroup,  the cycle in $\sum_S$ consisting of $\{x_0=0,x_1=0,x_2=0\}$ satisfies the assumption of (\ref{criterion}). So $(\P^2,G)$ is not birational to any action with IFP.
\end{proof}

 Dolgachev and Iskovskikh classify all such finite $G$ (\cite{dol}, Lemma 4.6). To find all transitive actions which are birational to the ones with IFP, we need to find all $G$ which do not contain any noncyclic subgroup. They are listed as following, which gives subtable of (\cite{dol}, Lemma 4.5 and Form 4.1). Here we denote $\tilde{G}$ to be the preimage of $G$ in $\C^*\times \SL_2(\C)$.

\begin{enumerate}
\item $G \cong (\Z/mk,\Z/m,\Z/nk,\Z/n)_{\alpha}$ ($\alpha$ is an automorphism of $\Z/k$), $\gcd(m,n)=1$;
\item $\tilde{G}\cong (\mu_{2m},\mu_{2m},H,H),G\cong \mu_m\times H$, where $H$ is a nonabelian binary polyhedral group,  $\gcd(m,|H|)=1$;
\item $\tilde{G}\cong (\mu_{6m},\mu_{2m},\bar{T},Q_8),G\cong (\mu_m\times Q_8).\Z/3$, $\gcd(m,2)=1$;
\item $\tilde{G}\cong (\mu_{4m},\mu_{2m}, Q_{4n},\Z/n), G\cong (\mu_{2m}\times \mu_{n}).\Z/2$. $m$ is even, $\gcd(m,n)=1$;
\item $\tilde{G}\cong (\mu_{4m},\mu_{m},Q_{4n},\Z/n),G\cong \mu_m\times D_{2n}$, $n$ is odd, $\gcd(m,n)=1$;
\item $\tilde{G}\cong (\mu_{4m},\mu_{2m},Q_{8n},Q_{4n}), G\cong (\mu_m\times Q_{4n}).2$, $\gcd(m,2n)=1$.\\
\end{enumerate} 

{\it Transitive imprimitive actions:}
\begin{lemma}[\cite{dol}, Theorem 4.7] Let $G$ be a transitive imprimitive finite subgroup of $\PGL_3$. Then
$G$ is conjugate to one of the following groups:
\begin{enumerate}
\item[$\bullet$] $G \cong  (\Z/n)^2 : Z/3$ generated by transformations
$$[\epsilon_nx_0, x_1, x_2], [x_0, \epsilon_nx_1, x_2], [x_2, x_0, x_1];$$
\item[$\bullet$] $G \cong (\Z/ n)^2 : S_3$ generated by transformations
$$[\epsilon_nx_0, x_1, x_2], [x_0, \epsilon_nx_1, x_2], [x_0, x_2, x_1], [x_2, x_0, x_1];$$
\item[$\bullet$] $G = G_{n,k,s} \cong (\Z/n \times \Z/\frac{n}{k}):\Z/3$, where $k > 1, k|n$ and $s^2 -s + 1 = 0$ mod $k$. It
is generated by transformations
$$[\epsilon_{n/k}x_0, x_1, x_2], [ \epsilon_n^sx_0, \epsilon_nx_1, x_2], [x_2, x_0, x_1];$$
\item[$\bullet$] $G \cong (\Z/n \times \Z/\frac{n}{3}):S_3$ generated by transformations
$$[\epsilon_{n/3}x_0, x_1, x_2], [\epsilon_n^2x_0, \epsilon_nx_1, x_2], [x_0, x_2, x_1], [x_1, x_0, x_2].$$
\end{enumerate}

\begin{ex}[$G=(\Z/3)^2:\Z/2)$]\label{example} $G$ is generated by 
$$[\epsilon_3^2x_0,\epsilon_3x_1,x_2], [x_2,x_0,x_1], [x_0,x_2,x_1].$$
$\sum_{\P^2}$ contains $9$ lines $\{x_i=\epsilon_3^k x_j\}$, which is the {\it Hessian arrangement}: each line passes through exact $4$ points of $\{ (1,\epsilon_3^i,\epsilon_3^j)(0\le i,j\le 2), (1,0,0),(0,1,0),(0,0,1)\} $. Through each of these $12$ points, there are exact $3$ lines in $\sum_{\P^2}$. Blowing up the above $12$ points, and then contracting the birational transform of $\sum_{\P^2}$, we get an action with IFP, and $3K_S$ is a trivial Cartier divisor. 
\end{ex}
 
\end{lemma}
\begin{prop}\label{transimpri}
All transitive imprimitive actions $(\P^2,G)$ which are birational to the ones with {\rm IFP} have $G$ as one of the following groups: $S_3$, $\Z/3:S_3\cong (\Z/3)^2:2$ and $\Z/n: \Z/3$.
\end{prop}
\begin{proof} 

When $G=(\Z/n)^2:K$ ($K=\Z/3$ or $S_3$) and $n> 1$, the subgroup action $(\P^2, (\Z/n)^2)$ is not birational to an action with only IFP. In fact, the cycle $\{x_0=0, x_1=0$ and $x_2=0\}\subset \sum_{\P^2}$ satisfies the assumption of (\ref{criterion}). When $n=1$, after possibly blowing up $(1,1,1)$, we can see the group action satisfies all the conditions of (\ref{cycstab}).

A similar argument shows when $G=G_{n,k,s}:\Z/3$, the action is not birational to an action with IFP if  $k\ne n$. For $\Z/n:\Z/3$, it acts on $\P^2$ with IFP.

For the last case, we only need to consider when $n=3$. And for this case $G=\Z/3:S_3$, we have an equivariant birational model with IFP as in (\ref{example}).\\
\end{proof}

{\it Primitive actions:} For the classical cases of finite primitive actions on $\P^2$ (a table of all such actions is given in \cite{Bl} or \cite{dol}, Theorem 4.8), we have the following result ,
\begin{prop}
\begin{enumerate}\label{primitive}
\item The action of the icosahedron group $A_5$ on $\P_2$ which leaves a nonsingular conic invariant is not birational to any action with {\rm IFP}.
\item The action of the Hessian group $\Hes\cong (\Z/3)^2:\bar{T}$ which is the automorphism group of the Hessian pencil 
$$x^3 + y^3 + z^3 + txyz = 0$$
is not birational to any action with {\rm IFP}.
\item The actions of the subgroups $G$ of the Hessian group, where $G=(\Z/3)^2:(\Z/4)$ and $(\Z/3)^2:Q_8$, are birational to actions with {\rm IFP}.
\item The action of the Klein group $L_2(7)$ of order 168 which is the automorphism group of the Klein quartic
$$x^3y + y^3z + z^3x=0 $$ 
is not birational to any action with {\rm IFP}.
\item The action of the Valentiner group of order 360 $(\cong A_6)$ which can be realized as
the full group of automorphisms of the nonsingular plane sextic
$$10x^3y^3 + 9zx^5 + 9zy^5 - 45x^2y^2z^2 - 135xyz^4 + 27z^6 = 0 $$ 
is not birational to any action with {\rm IFP}.
\end{enumerate}
\end{prop} 
\begin{proof}
We check the claim case by case:
\begin{enumerate} 
\item  $A_5$ acts on $\P^1$, so it acts on the complete linear system of $\sO(2)$ which is isomorphic to $\P^2$. Choose $(x_0^2,x_0x_1,x_1^2)$ to be the basis of $\P^2$. Now $A_5$ has a subgroup $D_4$ , whose nontrivial elements act on $\P^1$ as $(x_0,x_1)\to (x_1,x_0), (x_0,x_1)\to (x_1,-x_0)$ and $(x_0,x_1)\to (-x_0,x_1)$. Then the induced actions on $\P^2$ fix three lines $(y_0-y_2=0)$, $(y_0+y_2=0)$ and $(y_1=0)$ respectively. By (\ref{criterion}), we see that this action of $A_5$ is not birational to one with IFP.
\item The Hessian group $G_{216}$ has a homomorphism to $\PSL_2$ which induces the following exact sequence
$$1\to \Z_3:S_3\to G_{216}\to A_4 \to 1.$$ 
We can write a generator of $G_{216}$ as following: the kernel $\Z/3:S_3$ is the group as in (\ref{example}). And we have another 3 generators (cf. \cite{dol2}, 3.1.4):\\ 
$\sigma_1=\left(\begin{array}{ccc}
       1& 1 & 1\\
       1& \epsilon &\epsilon^2\\
     1&    \epsilon^2 & \epsilon  
\end{array}\right)$,
$\sigma_2=\left(\begin{array}{ccc}
       1& \epsilon  & \epsilon\\
       \epsilon^2& \epsilon &\epsilon^2\\
     \epsilon^2&    \epsilon^2 & \epsilon  
\end{array}\right)
$ and $\sigma_3=[\epsilon_3 x_0,x_2,x_1]$.

We notice that there is a subgroup of the Hessian group generated by $u: (x_0,x_1,x_2)\to(\omega x_0,x_1,x_2)$ and $v: (x_0,x_1,x_2)\to(x_0,\omega x_1,x_2)$. We know this $(\Z/3)^2$ action is not birational to an action with IFP according to (\ref{transimpri}).  
\item In the above exact sequence, $A_4$ has a subgroup $\Z/2$ which is generated by the image of $\sigma_1$ and a subgroup $(\Z/2)^2$ which is generated by the image of $\sigma_1$ and $\sigma_2$. Since $(\Z/3)^2:(\Z/4)$ is a subgroup of $(\Z/3)^2:Q_8$, we only need to prove the statement for the second case. Now the subgroup $(\Z/3)^2$ generated by $(x_0,x_1,x_2)\to(\omega x_0,\omega^2x_1,x_2)$, $(x_0,x_1,x_2)\to(x_2,x_0,x_1)$ is the only noncyclic abelian subgroup of $(\Z/3)^2:Q_8$. In this case, we can easily check that $(\Z/3)^2:Q_8$ acts on the model we construct in (\ref{example}), and gives an action with IFP.
\item By (\cite{dol2}, 6.5.2), the Klein quartic $x_0^3x_1+x_1^3x_2+x_2^3x_0=0$ is a specialization of the quartic 
$$C_{a,b,c}:x_0^4+x_1^4+x_2^4+ax_0^2x_1^2+bx_1^2x_2^2+cx_0^2x_2^2=0$$  
when $(a=b=c=\frac{-1+\sqrt{7}}{2})$. Notice that  $(\Z/2)^2$ acts on a general $C_{a,b,c}$, and the induced action on $\P^2$ is not birational to any action with IFP. Thus, the Klein group action is not birational to any action with IFP.
\item  By (\cite{Bl}), the Valentiner group will contain the icosahedron group as a subgroup, then by (1), we know it is not birational to any action with IFP.     
\end{enumerate}
\end{proof}

\noindent
{\bf The Case $\mathbf{S=\P^1\times \P^1}$}

$\Aut(S)=\PGL_2\wr S_2$. So every finite group of $\Aut(S)$ has a subgroup $G^0$ in $\PGL_2\times \PGL_2$ with the index at most 2. Let $G^0 \cap (\PGL_2\times \{e\})=G_1$, $G^0 \cap (\{e\}\times \PGL_2)=G_1$ and $H^0:=G_1\times G_2$.

\begin{prop}\label{F_0}
 If $G=G^0$, $(\F_0,G)$ is birational to an action with {\rm IFP} if and only if $|G_1|$ and $|G_2|$ are coprime.
\end{prop}

\begin{proof}
In the first case, if $G^0$ contains a subgroup of the form $G_1\times G_2$, and their cardinalities are not coprime, then it contains a subgroup conjugate to $\mu_n \times \mu_n$. For $\mu_n \times \mu_n$ acting on $\P^1 \times \P^1$, we observe that $\{x_0=0,x_1=0,y_0=0,y_1=0\}$ gives a cycle in $\sum_{\P^1 \times \P^1,G}$ which satisfies the assumptions of (\ref{criterion}). Conversely, it suffices to verify the three conditions of (\ref{cycstab}). Since a nonunit element $g\in G$ which fixes a curve must be in $G_1$ or $G_2$,  the curves in $\sum$ are fibers of one of the projections. The assumption $|G_1|$ and $|G_2|$ are coprime indeed implies that one of them, say $G_1$, is cyclic. Then there are at most two fibers of the form $pt \times \P^1$ in $\sum$. Then a curve $C\in \sum$ of the form $\P^1\times pt$ belongs to $\tilde{\sum}$.  We claim that any abelian subgroup of the stabilizer of a point is cyclic. In fact, such a stabilizer group will be isomorphic to a group of the form $(\Z/mk,\Z/m,\Z/nk,\Z/n)_{\alpha}$. Then the requirement of its abelian subgroups being cyclic is equivalent to  $\gcd(m,n)=1$, which is equivalent to the above coprimeness assumption. Hence, we can apply (\ref{cycstab}) to this case.
\end{proof}

\noindent We list all possible actions in the case when $G=G^0$,
\begin{enumerate}
\item $G_1\times G_2$, where $|G_1|$ and $|G_2|$ are coprime;
\item $(G,1,G,1)_{\alpha}=\{(g,\alpha(g))| G=\Z/n, D_{2n}, T, O$ or $I\}$;
\item $(D_{2m},\Z/m,O,T), \gcd(m,6)=1$;
\item $(D_{2m},\Z/m,O,(\Z/2)^2)$, $\gcd(m, 2) = 1$;
\item $(\Z/2m,\Z/m,O,T), \gcd(m,6)=1$;
\item $(\Z/3m,\Z/m,T,D_4), \gcd(m,2)=1$;
\item $(D_{2m},\Z/m,D_{4n},D_{2n}), \gcd(m,2n)=1$;
\item $(\Z_{2m},\Z/m,D_{2n},\Z_{n}), \gcd(m,n)=1$;
\item $(D_{2mk},\Z/m,D_{2nk},\Z/n)_{\alpha}, \gcd(m,n)=1$;
\item $(\Z/2m,\Z/m,D_{2n},D_{n}), \gcd(m,2n)=1$;
\item $(\Z/mk,\Z/m,\Z/nk,\Z/n)_{\alpha}, \gcd(m,n)=1$;
\end{enumerate}

For the argument later, we point out that in the case (9), when $k>2$ the group $G$ is isomorphic to $D_{2mnk}$; when $k=2$, there are 2 groups: besides $D_{4mn}$, there is another action which is birational to the action $(D_{2m},\Z/m,D_{4n},D_{2n})$ as in (7).  

\begin{prop}\label{f_0}
If $[G:G^0]=2$, $(\F_0,G)$ is birational to an action with {\rm IFP}, then $G$ is one of the following group: $\Z/2n$, $D_{2n}$($n$ is odd or $4k$), $F_{4n} (n=4k+2)$ or $G_{4n}(n=4k+2)$. When 
$(p^2\equiv -1 \mod n)$
 has a solution, we also have $H_{4n}$ , $I_{4n}$ ($n$ even), $J_{4n}$ ($n$ odd) where
\begin{enumerate}
\item $F_{4n}$ is the group generated by $(x,y)\to (-\frac{1}{x},-\frac{1}{y})$, $(x,y)\to (e^{\frac{2i\pi}{n}}x,e^{\frac{2i\pi}{n}}y)$ and $(x,y)\to (e^{\frac{i\pi}{n}}y,-e^{\frac{i\pi}{n}}x)$, and
\item $G_{4n}$ is the group generated by $(x,y)\to (-\frac{1}{x},-\frac{1}{y})$, $(x,y)\to (e^{\frac{2i\pi}{n}}x,e^{\frac{2i\pi}{n}}y)$  and $(x,y)\to (-y,x)$. 
\item $H_{4n}$ is the group generated by $(x,y)\to (-\frac{1}{x},-\frac{1}{y})$, $(x,y)\to (e^{\frac{2i\pi}{n}}x,e^{\frac{2pi\pi}{n}}y)$  and $(x,y)\to (-\frac{1}{y},x)$.
\item $I_{4n}$ is the group generated by $(x,y)\to (-\frac{1}{x},-\frac{1}{y})$, $(x,y)\to (e^{\frac{2i\pi}{n}}x,e^{\frac{2pi\pi}{n}}y)$  and $(x,y)\to (-\frac{1}{y},-x)$.
\item $J_{4n}$ is the group generated by $(x,y)\to (-\frac{1}{x},-\frac{1}{y})$, $(x,y)\to (e^{\frac{2i\pi}{n}}x,e^{\frac{2pi\pi}{n}}y)$  and $(x,y)\to (\frac{1}{y},-x)$.
\end{enumerate}

\end{prop}
\begin{proof}

In general $G^0$ has the form $(G,H, G,H)_{\alpha}$, where $G$ is given by the projection of $G^0$ on each factor.
As in the argument of (\ref{F_0}), we know the only possible case is when $H = 1$,
so $G = G^0\bullet\Z/2$, where $G^0$ is a polyhedral group. An
element $h \in G$ whose image in $\Z/2$ is nontrivial can be represented as $(h_1, h_2)\tau$,
where $\tau$ is the element of switching 2 factors. $(h_1, h_2)\tau$ fixes a curve if and only if $h_1=h_2^{-1}$, which is also equivalent to saying that $(h_1, h_2)\tau$ has order 2. In this case, it fixes the curve $(h_1(y),y)$. 
So if 
$$\mbox{($\ast$) $G$ contains a subgroup with the property: $H\cong (\Z/2)^2$ and $H\subsetneq G^0$},$$ 
then applying (\ref{criterion}), we know it is not birationally with IFP.

We claim if $G$ acts on $\F_0$ birationally with IFP, then $G^0$ is either cyclic or dihedral. Since both $I$ and $O$ contain $T$, We only need to rule out the case $G^0 = T$.
$\Aut(A_4)\cong S_4$, so after taking a conjugation of an element in $\PGL_2$, we can assume $\alpha=id$.  Now if $(g,g)\in G^0$, $h^{-1}(g,g)h = (h^{-1}_2 gh_2, h^{-1}_1 gh_1) \in G^0$,
so $h_1 h^{-1}_2$ is a commutator of $T$ in $\PGL_2$, which implies $h_1=h_2$. Furthermore, $h_1$ is in the
normalizer. So $G$ is either $T \times \Z/2$ $(h_1 \in T)$ or $O$. For both two cases, $G$
has a subgroup $H$ satisfying ($\ast$). 

$G^0 = (\Z/n, 1,\Z/n, 1)_s$, $\gcd(n, s) = 1$. Since $G$ is an extension of $\Z/2$ by
$\Z/n$, $G = \Z/2n$, $\Z/n + \Z/2$ ($n$ even) or $D_{2n}$. $\Z/n+\Z/2$ ($n$ even) and $D_{2n}$ ($n$ even) satisfying $(\ast)$. On the other hand, $\sum_{\Z/2n}$ is empty or a single curve. $\sum_{D_{2n}}$ ($n$ odd) consists of precise $n$ rational curves. Any two of them intersect at two indentical points. Blowing up these two points, we have a model satisfying (\ref{cycstab}). Hence, we conclude the action of $G$ on $\F_0$ is birationally with IFP if and only if $G = \Z/2n$ or $D_{4k+2}$.

If $G^0 = (D_{2n}, 1,D_{2n}, 1)_{\alpha}$ is generated by $a, b$ with $a^n = 1, b^2 = 1$. We will choose a representation: $a(z)=e^{\frac{2i\pi}{n}}z$ and $b(z)=-\frac{1}{z}$.  Then we can
assume $\alpha(b) = b$ by composing $\alpha$ with an action of an element in $\PGL_2$. When
$n = 2$, it is easy to see that $G = \Z/4+\Z/2$ or $Q_8$. So in the following argument,
we assume $n > 2$. $((h_1, h_2)\tau )^{-1}(a, \alpha(a))(h_1, h_2)\tau = (h^{-1}_2 \alpha(a)h_2, h_1^{-1}ah_1)$. Any conjugation of an element in $\PGL_2$ fixing $\Z/n$ will send $a$ to $a$ or $a^{-1}$. Then we can see $\alpha^4(a)=a$. If $\alpha(a) =a^{-1}$, after composing the conjugation of $b$, we can reduce to case that $\alpha=id$.

$\bullet \alpha=id:$

From now on, we change our notation by writing $h_1h_2^{-1}$ as $t$ and $h_2$ as $r$.  Now the above informations are read as: $t$ commutes with $D_{2n}$, $r$ normalizes $D_{2n}$ and $r^2 t\in D_{2n}$. Now if $t=e$, then we have a composition of group homomorphism, $f:G \to \PGL_2\times \Z/2 \to \PGL_2$. If $\Ker(f)=\Z/2$, then $\tau \in G$, and $G=D_{2n}\times \Z/2$. Otherwise, $f$ is an isomorphism from $G$ to its image, which is a polyhedral group containing $D_{2n}$ as an index $2$ subgroup.   Because of $(\ast)$, we have $(\F_0,G)$ is birational to an action with IFP if and only if $(G,G^0)$ is $(D_{4n}, D_{2n})$ and $n$ is even.

If $t\ne e$, if $r\in D_{2n}$, we can write $h$ to be $(-1,1)\tau$. $h^2=(-1,-1)$ implies $n$ is even. So $G=G_{4n}$. When $n$ is divided by 4, $(-i,i)\tau$ and $(i,i)\tau$ generate a subgroup as $H$ in $(\ast)$. On the other hand, $\sum_{G_{16k+8}}$ is empty. So in this case, it is birationally with IFP iff $n=4k+2$.  If $r\not\in D_{2n}$, since $-r^2\in G^0$. We can assume $r$ commutes with $a$, then $h^{-1}(b,b)h=(r^2b,r^2b)$. This implies $r^2 \in D_{2n}$. So $n$ is even. And the same argument as the previous case shows, $G$ is birationally with IFP if and only if $n=4k+2$.

$\bullet \alpha$ has order 4:

Now we can assume $h_1=bt_1$ and $h_2=t_2$, where $t_i$ ($i=1,2$) are commutators of $a$. So $((h_1,h_2)\tau)^2=(bt_1t_2,bt_2^{-1}t_1),$ which is never trivial. This is saying that $G$ acts on $F_0$ with IFP.
Assume $\alpha(a)=a^q$, then $q^2\equiv -1$ (mod $n$). Let $t_1t_2=a^{k_1}, t_2^2=a^{k_2}$, then $qk_1=k_1-k_2$. 

$n$ is even: it implies that $t_2$ is in $D_{2n}$ and  we can choose $t_2=e$. $\gcd(q-1,n)=2$ implies $k_1=0$ or $\frac{n}{2}$.  So we have 
\begin{enumerate}
 \item $H_{4n}$ generated by $a$, $b$ and $(b,1)\tau$.
\item $I_{4n}$ generated by $a, b$ and $(b,-1)\tau$.
\end{enumerate}
($H_{4n}\cong I_{4n}$ as abstract groups).

$n$ is odd: if $t_2$ is in $D_{2n}$, again we assume it is $e$. Since $\gcd(q-1,n)=1$, we know $t_1=e$, and we get $H_{4n}$. If $t_2$ is not in $D_{2n}$. We can assume it is $-1$, which implies that $t_1=-1$. Then we have a group $J_{4n}$ generated by $a$, $b$ and $(-b,-1)\tau$. $J_{4n}$ is also isomorphic to $H_{4n}$ as an abstract group.
\end{proof}

\noindent
{\bf The Case S=$\F_e$, $e>1$}
\begin{prop}\label{F_n}
The actions on $\F_e (e\ge 1$), which are birationally with {\rm IFP}, are as follows:\\
{\rm (I)} When $e$ is even, $G$ are the groups in (\ref{intransitive});\\ 
{\rm (II)} When $e$ is odd, $G\subset \C^*\times \PGL_2$ are the groups on the above list with the form $(\Z/mk,\Z/m,G_1,G_2)$.  
\end{prop}
\begin{proof} 
 
In fact, $\Aut(\F_n) \cong \Aut(\P(1, 1, n)) \cong \C^{n+1} : GL_2/\mu_n$, so if a finite
group $G \in Aut(\F_n)$, then $G \in GL_2/\mu_n$. Since
$$1 \to \Z/2 \to \C^* \times \SL_2 \to \GL_2 \to 1,$$
we have $GL_2/\mu_n \cong C^*/\mu_n\times\PGL_2 \cong \C^*\times \PGL_2$ ($n$ is even) or $(\C^*/\mu_n \times \SL_2)/\Z/2 \cong \GL_2$ ($n$ is odd). Use an argument similar to the above proof,
we have the conclusion. The details are left to the reader.
\end{proof}

\noindent
{\bf The Case $S$ is Nonminimal $G$-Ruled Surface}

\begin{prop}
Assume a minimal action $(S,G)$ satisfies that $\pi:S\to P^1$ gives a $G$-equivariant fiberation, and $S$ itself is not minimal. Then it is not birational to any action with {\rm IFP}.
\end{prop}
\begin{proof}
Consider the natural group homomorphism $f: G \to \Aut(\Pic(S))$. 

If its kernel $G_0$ is not trivial, then by (\cite{dol}, Proposition 5.5), we know $S$ is an {\it exceptional conic bundle}, i.e. the minimal resolution of the degree $2g+2$ hypersurface 
$$ F_{2g+2}(T_0,T_1)+T_2T_3=0$$ 
in the weighted projective space $\P(1,1,g+1,g+1)$ (See the construction of Section 5.2 of \cite{dol}). 
The automorphism group of $S$ is an extension of $P$ by $N$, where $P$ is the subgroup of $\PGL_2$ leaving the zero sets of $F_{2g+2}$ invariant and $N\cong \C^*:2$ is a group of matrices with determinant $\pm 1$ leaving $T_2T_3$ invariant (cf. \cite{dol} Proposition 5.3). The kernel $K$ of $f$ is the intersection of $G$ with $N$ and fixes the coordinates $T_0$, $T_1$ and left $T_2T_3$ invariant. So it fixes the curve which is the birational transform of $C: F_{2g+2}(T_0,T_1)+T^2=0$ pointwise. Since this curve is of genus greater or equal to 1, the action is not birationally with {\rm IFP}.

If $G_0=\{ e \}$. Then thanks to the following proposition, we know $G$ contains an order-$2$ element fixing a curve of genus $g\ge 1$ pointwise.

\begin{prop}[\cite{dol}, Theorem 5.7]
Let $G$ be a minimal finite group of automorphisms of a conic bundle
$\varphi : S \to \P^1$ with a set $\Sigma$ of singular fibres. Assume $G_0 = {e}$. Then $k = |\Sigma| > 2$
and one of the following cases occurs:
\begin{enumerate}
 \item 
 $G = 2\bullet P$, where the central involution $h$ fixes pointwise an irreducible smooth
bisection $C$ of $\varphi$ and switches the components in all fibres. The curve $C$ is
a curve of genus $g = (k − 2)/2$. The conic bundle projection defines a $g^1_2$
on C with ramification points equal to singular points of fibres. The group
$P$ is isomorphic to the group of automorphisms of $C$ modulo the involution
defined by the $g^1_2$.
\item (2) $G \cong 2^2 \bullet P$, each nontrivial element $g_i$ of the subgroup $2^2$ fixes pointwise
an irreducible smooth bisection $C_i$. The set $\Sigma$ is partitioned in 3 subsets
$\Sigma_1,\Sigma_2,\Sigma_3$ such that the projection $φ\varphi : C_i → \P^1$ ramifies over $\Sigma_j + \Sigma_k$, $
i \neq j \neq k$. The group $P$ is subgroup of $Aut(\P^1)$ leaving the set $\Sigma$ and its
partition into 3 subsets $\Sigma_i$ invariant.
 \end{enumerate}
\end{prop}
\end{proof}

It remains to study the cases when $(S,G)$ is minimal and $S$ is a smooth del Pezzo Surface, For the Del Pezzo surface of degree 7 and 8, there does not exsit any minimal action.

\noindent {\bf The Case $S$ is the del Pezzo Surface of Degree 6}
\begin{prop}[\cite{dol}, Theorem 6.3]
Let $G$ be a minimal subgroup of a del Pezzo surface $S$ of degree 6. Then
$G = H\bullet \langle s \rangle$,
where $H$ is an transitive imprimitive finite subgroup of $\PGL_3$ and $s$ is 
the lift of the standard quadratic transformation.
\end{prop}

\begin{prop}\label{degree6}
Notation as above. If $(S,G)$ is birational to an action with {\rm IFP}, then  $G$ is $S_3$ or $\Z/2\times (\Z/n:\Z/3)$ for $n \ge 1$.
\end{prop}
\begin{proof}
The $G=H\bullet \langle s \rangle$ acting on $S$ birationally with IFP implies the same thing holds for $H$ on $S$, which has a minimal model of $H$ on $\P^2$ as an imprimitive action. So $H$ can be only the groups as in (\ref{transimpri}).

We claim the action of $G=\Z/2\times S_3$ ($H=S_3$) is not birational to any action with IFP. In fact, the abelian subgroup generated by the lifting of the order $2$ element $\rho:(x_0,x_1,x_2)\to (x_1,x_0,x_2)$ and the Cremona transformation $\tau$ is isomorphic to $(\Z/2)^2$. We have $S^{\rho}=(x_0=x_1)$, $S^{\rho \tau}= (x_0x_1=x_2^2)$. This gives a cycle in $\sum_{S} $ satisfying the assumption of (\ref{criterion}). 

The only remaining case is when $H=\Z/n:\Z/3$. We claim the action $(S,G)$ itself is already with IFP. First, we know any nontivial element in the subgroup $G_{n,s}$ acts on $S$ with IFP. For an element $g\in G_{n,s}:\Z/6$, $g^6 \in G_{n,s}$, if $g^6$ is nontrivial, we have $g$ acts on $S$ with IFP. $n| s^2-s+1\equiv 0$ implies $n$ is odd. Hence, we know if $S^g$ is not isolated, we have $g^3$ is trivial. Thus, we only need to verify the statement for the case when $G=G_{3,2}:\Z/6=(\Z/3)^2 \times \Z/2$. However, both the Cremona transformation and any element in $(\Z/3)^2$ fix finite points, so we conclude that the $G_{n,s}:\Z/6$ action is an action with IFP.

\end{proof}   

\noindent{\bf The Case $S$ is the del Pezzo Surface of Degree 5}

$\Aut(S)=S_5$, we assume that we get $S$ by blowing up 4 points $(1,0,0), (0,1,0), (0,0,1)$ and $(1,1,1)$. 

\begin{prop}[\cite{dol}, Theorem 6.4]  Let $(S,G)$ be a minimal Del Pezzo surface of degree $d = 5$. Then
$G = S_5, A_5, \Z/5 : \Z/4, \Z/5 : \Z/2$, or $\Z/5$.
 \end{prop}

\begin{prop}\label{degree5}
Let $(S,G)$ be a minimal action on the smooth Del Pezzo surface $S$ of degree 5, and assume $(S,G)$ is birational to an action with {\rm IFP}. Then $G=\Z/5:\Z/4, \Z/5:\Z/2$ or $\Z/5$.
\end{prop}
\begin{proof}
For $G=A_5$: there is a cycle in $\sum_S$ whose edges are birational transforms of lines $(x_i+x_j=x_k)(\{i,j,k\}=\{1,2,3\})$. The vertices are $(1,1,0),(1,0,1),(0,1,1)$ and each of them has a noncyclic abelian stabilizer. Then by (\ref{criterion}), we know that $A_5$ is not birational to any action with IFP. 

$G=\Z/5:\Z/4$:  We can represent the elements of $G \cong \Z/5:\Z/4=D_{10}:\Z/2$ as follows:
\begin{enumerate}
\item[$\bullet$] $(12345)=(x_0,x_1,x_2)\to(x_0(x_2-x_1),x_2(x_0-x_1),x_0x_2)$;
\item[$\bullet$] $(2354)=(x_0,x_1,x_2)\to(x_2(x_0-x_1),x_2(x_0-x_2),x_1(x_0-x_2))$.
\end{enumerate}
Then $\sum_S$ contains $5$ irreducible curves, any $2$ of which intersect at $2$ identical points $(\frac{-1- \sqrt5}{2},\frac{3+ \sqrt5}{2},1)$ and $(\frac{-1+ \sqrt5}{2},\frac{3- \sqrt5}{2},1)$. We can first blow up these $2$ points, then contract the birational transforms of the above $5$ curves, then it gives an action on $\F_0$.

\end{proof}

In the last part of this section, we will prove for any minimal action on a smooth del Pezzo surface of degree less than 5, the group always contains a nontrivial element which fixes a curve of genus $g>0$ pointwise. In particular, it implies there does not exist any minimal action on such surfaces, which is birational to an action with IFP. First, we study some general theory of a finite group $G$ acting on a smooth surface $S$, and apply it to the case when $S$ is rational.
 
For any nontrivial automorphism $g$ of a surface $S$, by the Lefschetz fixed-point formula, we have 
$$2-\Tr_1(g)+\Tr_2(g)=e(S^g)=s+\sum_{j=1}^{t}(2-2g(C_j)).$$
Here $\Tr_i$ means the trace of the $g$ acting on the $i$-th singular cohomology. $C_j$ are the disjoint smooth curves fixed by $g$, and $s$ is the number of the isolated fixed points. When $S$ is rational, $H^1(S)=0$, and $\Pic(S)=H^2(S,\Z)$. Then for a group $G$ acting on $S$, we have

{\begin{eqnarray*}
\rank (\Pic(S)^G)& =& \frac{1}{|G|}\sum_{g\in G}\Tr_2(g)\\
                 & =& \frac{1}{|G|}(\rank(S)+\sum_{g\in G-\{e\}}(s-2+\sum_{j=1}^{t}(2-2g(C_j))),
\end{eqnarray*}

\noindent If $(S,G)$ is minimal, but $S$ is not a $G$-equivariant conic bundle, then we have $\rank (\Pic(S)^G)=1$.

\begin{rem}
Although the above general theory is illuminating, in the following proofs, we have to use the classification results from \cite{dol}. For this reason, we will use the terminology there without referring. However, it would be nice to find a straightforward argument which does not heavily depend on the classification results. 
\end{rem}

\noindent {\bf The Case $S$ is a del Pezzo Surface of Degree 4}

The reader can check Subsection (6.4) of \cite{dol} for all minimal actions $G$ on $S$, which is a del Pezzo surface of degree 4. Let us summarize the facts we need here. $S$ is isomorphic to a nonsingular surface of degree 4 in $\P^4$ given by equations
$$ F_1 =\sum_{i=0}^4 T^2_i = 0, F_2 =\sum^4_{i=0} a_iT^2_i = 0,$$ 
where all $a_i$’s are distinct. The Weyl group is $W(D_5)=\Z/2^4:S_5$ and the automorphism group of $S$, which is a subgroup of $W(D_5)$, always has $(\Z/2)^4$ as a normal subgroup.
It is given by changing even number of signs of the coordinates. So for $A$ an even cardinality subset of $\{0,1,2,3,4\}$, we can compute the fixed locus of $i_A$, where $i_A$ is the automorphism of changing the signs of coordinates corresponding to $A$.  When $|A|=4$, the fixed locus is an elliptic curve, and when $|A|=2$, the fixed locus consists of isolated points.  $(\Z/2)^4\cap G$ can only be $e,i_{ab}$ or $\langle i_{ab},i_{ac} \rangle$. 

The subgroup $G'$ of $\Aut(S)$ can be realized as the stabilizer of a set of 5 skew lines on $S$. Thus $G'$ is isomorphic to a group of projective transformations of $\P^2$ leaving invariant a set of 5 points. Since there is a unique conic through these
points, the group is isomorphic to a finite group of $\PGL_2$ leaving invariant a set
of 5 distinct points. It follows that a subgroup leaves
invariant a set of 5 distinct points if and only if it is one of the following groups
$\Z/2,\Z/3,\Z/4,\Z/5, S_3,D_{10}$.

\begin{prop}[\cite{dol}, Theorem 6.9]
 Let $(S,G)$ be a minimal del Pezzo surface of degree 4. Then
$G$ is isomorphic to one of the following groups:
\begin{enumerate}
 \item $\Aut(S) \cong (\Z/2)^4$,
$$(\Z/2)^4, (\Z/2)^3, (\Z/2)^2.$$
\item $\Aut(S) \cong (\Z/2)^4 : \Z/2$,
$$\Z/2 \times \Z/4, D_8, L_{16}, (\Z/2)^4 : \Z/2,$$
and from the previous case.
\item $\Aut(S) \cong (\Z/2)^4 : \Z/4$,
$$8, (\Z/2)^2 : \Z/8, (\Z/2)^4 : \Z/4,$$
and from the previous two cases.
\item $\Aut(S) \cong (\Z/2)^4 : S_3$,
$$(\Z/2)^2 \times \Z/3, \Z/2 \times A_4(\cong (\Z/2)^3:\Z/3), (\Z/2)^4 : \Z/3, (\Z/2)^4 : S_3,$$
and from Cases 1) and 2).
\item  $\Aut(S) \cong (\Z/2)^4 : D_{10}$,
$$(\Z/2)^4 : D_{10}, (\Z/2)^4 : \Z/5,$$
and from Cases 1) and 2).
\end{enumerate}
\end{prop}

We claim for all these minimal actions, it contains an element of the form $i_{abcd}$. Otherwise,  $K:=(\Z/2)^4\cap G$ can only be $e,i_{ab}$ or $\langle i_{ab},i_{ac} \rangle$. Since all of these actions are not minimal, we conclude that $G\subsetneq (\Z/2)^4$, hence the image $G'$ of $G$ in $S_5$ is non-trivial. Also we need to only look at the cases that $|K|\le 4$. In Case (2): $G=\Z/2 \times \Z/4\Rightarrow K=\langle i_{01},i_{0123}\rangle$; $G=D_8 \Rightarrow K=\langle i_{0123},i_{0124}\rangle$. In Case (3): $G=\Z/2 \times \Z/4 \Rightarrow K=\langle i_{0123}\rangle$. In Case (4), $G=(\Z/2)^2 \times \Z/3 \Rightarrow K=\langle i_{0123},i_{0124}\rangle$.

Thus we conclude

\begin{prop}
Any minimal actions on a Del Pezzo surface $S$ of degree 4 is not birational to an action with {\rm IFP}.
\end{prop}

\noindent
{\bf The Case $S$ is a Cubic Surface}

All the possible minimal actions on a smooth cubic surface are classified in Subsection (6.5) of \cite{dol}. There are cyclic groups, whose generator fixes an elliptic curve pointwise:
\begin{enumerate}
 \item [$\bullet$] $[t_0, t_1, t_2, -t_3]$
$$F = T_3^2 L_1(T_0, T_1, T_2) + T^3_0 + T^3_1 + T^3_2 + \alpha T_0T_1T_2.$$ 
The generator is of type $4A_1$.
 \item [$\bullet$] $[t_0, t_1, t_2, \epsilon_3t_3]$
$$F = T^3_0 + T^3_1 + T^3_2 + T^3_3 + \alpha T_0T_1T_2.$$
 The generator is of type $3A_1$.
\end{enumerate}
We will check that for all minimal actions $(S,G)$, $G$ always contain elements of one of the above types.

Table 4 of \cite{dol} gives all types of cubic surfaces and their automorphic groups. Because there are some specializations as:
$$\mbox{ IV} \to \mbox{\rm III}, \mbox{\rm IV }\to \mbox{I}, \mbox{ } \mbox{ } \mbox{VI, VIII,IX} \to \mbox{I}, \mbox{ } \mbox{XI} \to \mbox{X},$$
it suffices to consider the surfaces of types {\rm I, II, III, V, VII, X}.

\begin{prop}[\cite{dol}, Theorem 6.14]
Let $G$ be a minimal subgroup of automorphisms of a nonsingular
cubic surface of type {\rm I, II, III, V, VII, X}. Then $G$ is isomorphic to one of the following groups (The number $n$ in $G(n)$ means there are $n$ different conjugacy classes.:
\begin{enumerate}
 \item $G$ is a subgroup of automorphisms of a surface of type {\rm I}.
$$S_4(3), S_3 (2), S_3 \times \Z/2, S_3 \times \Z/3 (2), (\Z/3)^2 : (\Z/2) (2), (\Z/3)^2 : (\Z/2)^2,$$
$$H_3(3) : \Z/2, H_3(3), (\Z/3)^3 : (\Z/2) (2), (\Z/3)^3 : (\Z/2)^2 (2),$$
$$(\Z/ 3)^3 :\Z/ 3,(\Z/3)^3 : S_3,(\Z/ 3)^3:D_8, (\Z/3)^3 : S_4, (\Z/3)^3 : \Z/4,$$
$$(\Z/3)^3, (\Z/3)^2 (3), (\Z/3)^2 \times \Z/2, \Z/9, \Z/6 (2), \Z/3.$$

\item  $G$ is a subgroup of automorphisms of a surface of type {\rm II}.
$$S_5, S_4.$$
\item $G$ is a subgroup of automorphisms of a surface of type {\rm III}.
$$H_3(3) : \Z/4, H_3(3) : \Z/2, H_3(3), S_3 \times \Z/3, S_3, (\Z/3)^2, \Z/12, \Z/6, \Z/3.$$
\item $G$ is a subgroup of automorphisms of a surface of type {\rm IV}.
$$H_3(3) : 2, H_3(3), S_3 (2), 3 \times S_3 (2), (\Z/3)^2 (2),\Z/ 6,\Z/ 3.$$
\item $G$ is a subgroup of automorphisms of a surface of type {\rm V}.
$$S_4, S_3.$$

\end{enumerate}
\end{prop}

We will discuss them case by case:

\begin{enumerate}

\item [(i)] The cubic surface of type \rm{I}: 
$$T_0^3+T_1^3+T_2^3+T_3^3=0.$$ 
Its automorphism group of is $(\Z/3)^3:S_4$. The factor of $(\Z/3)^3$ comes from sending $(T_0,T_1,T_2,T_3)$ to $(\epsilon_3^aT_0,\epsilon_3^bT_1, \epsilon_3^cT_2,\epsilon_3^bT_3)$, which we denote it as $[a,b,c,d]$. We require $$a+b+c+d=0 \mbox{ in } \F_3.$$ Then up to symmetry, there are 3 different type of 1-dimensional space in $\F_3^3$, which are represented by 
$$ [1,2,0,0], [1,1,1,0],[1,1,2,2], $$
The different classes of 2-dimenionsal spaces of $\F_3^3$ are given be the orthogonal complement in $(\F^3)^3$ with respect to the dot-product pairing on $\F_3^4$. Let $K=G\cap (\Z/3)^3$, and $\dim_{\F_3}K=k$. 
\begin{enumerate}
\item $k=0$, $G$ is either $S_4$ or $S_3$. When $G=S_4$, it has 3 different conjugacy classes. Each of them has the 6 elements of type $4A_1$; when $G=S_3$, it has 3 elements of type $4A_1$;

\item $k=1$, when $K=\langle [1,1,2,2]\rangle$, the action is not minimal; when $K=\langle [1,1,1,0]\rangle$,  the generator of $K$  is of type $3A_2$; and when $K=\langle [0,0,1,2]\rangle$, $G=S_3$ and it contains 3 elements of type $4A_1$;
\item $k=2$, if $K$ is an orthogonal complement of $[1,1,2,2]$, then the action is not minimal. But for other 2 cases,  $K$ contains element of type $[1,1,1,0]$, which is indeed of type $3A_1$;
\item $k=3$, $K$ contains a $3A_2$ type element.
\end{enumerate}
\item[(ii)] The cubic surface of type \rm{II}:
$$T_0^2T_1+T_1^2T_2+T_2^2T_3+T_3^2T_0=0.$$
The surface is isomorphic to the Clebsch diagonal cubic surface in $\P^4$ given by
the equations
$$\sum^4_{i=0}T^3_i=\sum^4_{i=0}T_i = 0.$$
The group $S_5$ acts by permuting the coordinates. The transposition $(12)$ is of type $4A_1$, and the minimal action, which has $G= S_5$ or $S_4$, always contains transpositions.

\item[(iii)] For the cubic surface of type \rm{IV}:
$$T^3_0 + T^3_1 + T^3_2 + T^3_3 + 6aT_1T_2T_3=0,$$ 
where the parameter $a$ satisfies $a − a^4 \neq 0$, $8a^3 \neq −1$, and $20a^3 + 8a^6 \neq 1$. Its automorphic group is $H_3(3) : 2$. When $a=1$, its specialization is of type \rm{I}. Given a smooth family of cubic $G$-surfaces $(S_t,G)$, both the minimality of the action and the type of an element $g\in G$ are invariant.  Hence from the argument in \rm{I}, all the minimal actions on type \rm{IV} cubic surface are not birational to actions with IFP.

\item[(iv)] For the cubic surface of type \rm{III}:
$$T^3_0 + T^3_1 + T^3_2 + T^3_3 + 6aT_1T_2T_3=0, $$
where $20a^3 + 8a^6 =1$. It is obviously a specialization of type \rm{IV}, and there are two new groups: $H_3(3):4$ and $\Z/12$. Both of them contain elements of type $3A_1$. 

\item[(v)] For the cubic surface of type \rm{V}, whose equation is 
$$T_0^3+T_0(T_1^2+T_2^2+T_3^2)+aT_1T_2T_3=0,$$
where $9a^3 \neq 8a, 8a^3 \neq −1$,
Its automorphic group is  $S_4 \cong (\Z/2)^2 : S_3$ acting by permuting the coordinates $T_1, T_2, T_3$ and multiplying them by $-1$ leaving the monomial $T_1T_2T_3$ unchanged.The only minimal action is when $G=S_4$ or it subgroup $S_3$. We notice that in the above representation, any transposition contained in $S_3$ is of type $4A_1$.  

\end{enumerate}

For a cubic surface of any other type, if it has some minimal actions, then it can be specialized to one of the above types. As in (\rm{iii}), we conclude none of them is birational to an action with fixed points. Thus we conclude,
\begin{prop}
Any minimal action on a cubic surface $S$ is not birational to an action with {\rm IFP}.
\end{prop}

\noindent {\bf The Case $S$ is a del Pezzo Surface of Degree 2}

All the minimal actions on a del Pezzo surface $S$ of degree 2 are discussed in Subsection (6.6) of \cite{dol}. Since $S$ can be written as a double curve of $\P^2$ branched over a quartic curve $B$. We know that there is a homomorphism $\Aut(S)\to \Aut(B)$ with the kernel generated by the Geiser involution. There are also cyclic groups, whose generator fixes an elliptic curve pointwise:
\begin{enumerate}
 \item [$\bullet$] $[t_0, t_1, -t_2, t_3]$
$$F = T_3^2+T_2^4 +T_2^2L_2(T_0, T_1) + L_4(T_0,T_1).$$ 
The generator is of type $4A_1$.
 \item [$\bullet$] $[t_0, t_1, \epsilon_3t_2,t_3]$
$$F = T_3^2 + T_2^3L_1(T_0,T_1) + L_4(T_0,T_1).$$
 The generator is of type $3A_1$.
\end{enumerate}
Then it suffices to check any minimal action $(S,G)$ satisifes that $G$ contains either the Geiser involution or an element of type $4A_1$ or $3A_1$. From \cite{dol} Lemma 6.16, we notice that if $g$ is an element of order 4, then $g^2$ is an element of type $4A_1$ or the Geiser involution, depending on the image of $g$ in $\Aut(B)$ is of order 4 or 2; if $g$ is an element of order 6, $g^2$ is of type $3A_1$ or $g^3$ is the Geiser involution, depending on the image of $g$ in $\Aut(B)$ is of order 6 or 3.  Thus a minimal group contain an element of order 4 or 6  is not birationally with IFP. 

\begin{prop}[\cite{dol}, Theorem 6.17] Let $G$ be a minimal group of automorphisms of a Del Pezzo surface
of degree 2. Then $G$ either contains the Geiser involution or equal to one of the minimal lifts of a subgroup of $\Aut(B)$ as the following:
\begin{enumerate}
 \item Type I: $L_2(7), S_4(2), D_8$

\item Type II: $\Z/4^2:S_3 (2), S_4 (2), (\Z/4)^2 :\Z/3, A_4 , (\Z/4)^2 :\Z/2(3),$ 
$$M_{16}, AS_{16} (2),D_8 ,(\Z/4_2),(\Z/2 \times \Z/4) (2), \Z/4$$

where $AS_{16}$ is the group of a presentation: 
$$a^4 = b^2 = c^2 = [a, b] = 1, [c, b]a^{-2} = [c, a] = 1.$$
$M_{16}$ is the group of a presentation: 
$$a^8 = b^2 = 1, [a, b]a^4 = 1.$$

\item Type III: 
$$\Z/4\bullet A_4 (2), D_8:3, AS_{16}(2),D_8, \Z/12, \Z/6,\Z/2 \times \Z/4, \Z/4.$$
\item Type IV: $S_4, D_8$.
\item Type V: $AS{16}(2),D_8, \Z/2 \times \Z/4 (2),\Z/4$.
\item Type VII: $D_8$.
\item Type VIII: $\Z/6$.

 \end{enumerate}
\end{prop}

Since there is a specialization: 
$$\mbox{IX} \to \mbox{IV} →\to \mbox{I, II}, \mbox{ XII} \to \mbox{X} \to \mbox{VII} \to\mbox{V} \to \mbox{II, III},
\mbox{ XI} →\to \mbox{VIII} \to \mbox{III}.$$
We only need to discuss the minimal groups which do not contain the Geiser involution and for surfaces of type I, II or III.
\begin{enumerate}
 \item [(i)] For the del Pezzo surface of degree 2 of type I:
$$T_3^2 + T_0^3 T_1 + T_1^3 T_2 + T_2^3 T_0=0.$$
The order-2 element in $L_2(7)$ is of type $4A_1$. $S_4$ and $D_8$ contain elements of order 4.
\item [(ii)] For the del Pezzo surface of degree 2 of type II:
$$T_3^2 + T_0^4 + T_1^4 + T_2^4=0.$$
All the groups contain an element of order 4.
\item [(iii)] For the del Pezzo surface of degree 2 of type III:
$$T_3^2+ T_2^4+ T_0^4+ aT_0^2 T_1^2 + T_1^4=0 (a^2=-12).$$
The groups contain an element of either order 4 or 6.
\end{enumerate}

Thus we conclude
\begin{prop}
Any minimal actions on a Del Pezzo surfaces $S$ of degree 2 is not birational to an action with {\rm IFP}.
\end{prop}

\noindent {\bf The Case $S$ is a del Pezzo Surface of Degree 1}

For the remaining cases of del Pezzo surfaces of degree 1. The idea is similar as the case of cubic surfaces. We will only schetch the proof and leave the details to the reader.

All the minimal actions on a del Pezzo surface $S$ of degree 1 are discussed in Subsection (6.7) of \cite{dol}. 
Any $S$ can be written as a degree $6$ hypersuface in the weighted projective space $\P(1,1,2,3)$ with the equation:  
$$S: T_3^2+T_2^3+T_2L_4(T_0,T_1)+L_6(T_0,T_1)=0.$$ 

The order-2 element of $G$ is one of the following:
\begin{enumerate}
 \item The Bertini involution:$[t_0,t_1,t_2,-t_3]$, which fixes a genus-4 curve pointwise;
\item $[it_0,-it_1,-t_2,it_3]$, where $F_4 = F_2(T_0^2, T_1^2) \neq 0, F_6 = F_3(T^0_2, T_1^2)$. It fixes the genus-1 curve $S\cap(t_0=0)$  pointwise; and 
\item  $[-t_1, t_0,-t_2, it_3]$, where
$F_4 = a(T_0^4+ T_1^4) + bT_0^2 T_1^2 , F_6 = a(T_0^6 − T_1^6) + bT_0T_1(T_0^4 + T_1^4)$. It fixes the genus-1 curve $S\cap (t_1= it_0).$

\end{enumerate}

Therefore, if the action is birationally with IFP, it cannot contain any order-$2$ elements, which means the order of the group is odd. We can also assume that $G$ does not contain the order-3 element $[t_0,t_1,\epsilon_3t_2,t_3]$ for it fixes a genus-2 curve pointwise. Then there are only a very small number of cases remaining:
\begin{enumerate}
 \item Type I, II, VII, XV: such group does not exsit; and
\item Type IV, VIII: $\Z/5$ generated by $[t_0,\epsilon_5t_1,t_2,t_3]$, which fixes a genus-1 curve pointwise. 

\end{enumerate}

 Thus we conclude
\begin{prop}
Any minimal actions on Del Pezzo surfaces $S$ of degree 1 is not birational to an action with {\rm IFP}.
\end{prop}

\section{group action on log del pezzo surface}
In this section, we aim to prove Theorem(\ref{main2}), namely given a finite group $G$, assume $G$ can act on a rational surface $\tilde{S}$ containing (at worst) quotient singularities which gives an action with only IFP, we would like to determine whether we can choose $\tilde{S}$ to be a log del Pezzo surface. The idea is to run the equivariant minimal model program for the pair $(S,G)$. For the general theory of minimal model program, see e.g. \cite{kollarmori}. It is well-known that a normal surface singularity is $klt$ if and only if it is a quotient singularity (cf. \cite{kollarmori}, 4.18). Thus if we start with a surface which contains (at worst) quotient singularities, and run the minimal model program, after a sequence of divisorial contractions, we still have a surface of the same type singularities. 

\begin{lemma}
If $(S,G)$ is an action with {\rm IFP}, and $R\cong S/G$, then running the equivariant minimal model program for $(S,G)$ is equivalent to running the (ordinary) minimal model program for $R$. 
\end{lemma}
\begin{proof}
The morphism $\pi: S\to R$ is finite, $\overline{NE}(S)^{G}=\overline{NE}(R)$ and $\Pic(S)^{G}=\Pic(R)$. So it suffices to prove $\pi^*(nK_{R})=nK_S$ for sme integer $n$ such that $nK_{R}$ is Cartier. In fact, after removing those isolated branched points, $\pi$ is an \'etale morphism, so the equality holds in this case. Hence, we can conclude $\pi^*(nK_{R})=nK_S$, since this is an equality of divisors.
\end{proof}

\begin{prop}\label{mmp}
Consider the groups $G$ which can act on a smooth rational surface $\tilde{S}$, such that it is birational to an action $(S,G)$ with {\rm IFP} and $K_S$ is not $\Q$-effective. These groups are precisely the groups which can act on some del Pezzo Surfaces with {\rm IFP}. In other words, they are precisely all quotient groups of $\pi_1$ of smooth loci of log Del Pezzo surfaces.
\end{prop}  

\begin{proof} For a surface with only quotient singularities, we know the $\Q$-effectivity of the canonical class is equivalent to the pseudo-effectivity. 
Let $R$ be a log del Pezzo surface, $G$ be a finite quotient group of $\pi_1(R^{sm})$, and $S$ be the corresponding cover over $R$, which is branched at finite points. $(S,G)$ gives us an action with IFP. As $S$ is also a log del Pezzo surface, $K_S$ is not pseudo-effective.

Conversely, we start with $(S,G)$ which is an action with IFP. Take $R \cong S/G$. We know $K_R$ is not pseudo-effective either. Running a log minimal model program for $R$, 
$$R=R_0\to R_1\to R_2\to \cdots \to R_n,$$
$R$ and $R_i$ are birational. The minimal model program preserves the non-effectivity assumption of $\Q$-divisor $K_R$. Thus it terminates with a Fano contraction to a lower dimensional variety. If it contracts to a point, which is equivalent to saying $\rho(R_n)=1$, then $R_n$ is a log del Pezzo surface, so the group $G$ is a quotient group of $\pi_1(R^{sm})$ which itself is a quotient group of $\pi_1(R_n^{sm})$ by (cf. \cite{keel}, 7.3). Otherwise, $R_n$ contracts to $\P^1$. By the above lemma, if we look at the corresponding $G$-equivariant minimal model program, 
$$S=S_0\to S_1\to S_2\to \cdots \to S_n,$$
it gives a contraction from $S_n$ to $\P^1$, which is a $G$-equivariant fiberation. Let $G_0$ be the kernel of the natural group homomorphism $\rho:G\to G_{|\P^1}$. For any element $g_0 \in G_0$, it acts on every fiber. In particular, for every fiber, the set of the fixed points is nonempty. Hence $g_0$ will fix some curves pointwise. Since $(S_n,G)$ is an action with IFP, we conclude that $G_0$ is trivial. So $G$ is a subgroup of $\PGL_2(\C)$, but any such group can diagonally act on $\P^1 \times \P^1$,  giving an action with IFP.
\end{proof}

The following criterion is useful to prove the non-effectivity.
\begin{lemma}\label{kodaira}
Let $S$ be a projective surface with a $G$ action which is birational to an action with {\rm IFP}. Suppose there exists a $G$-birational proper morphism $\phi: \bar{S}\to S$ with the property: for any $\Q$-divisor $E$, which supports on $\sum_{\bar{S}}$, we have $\kappa (K_{\bar{S}}+E)=-\infty$ provided $\rdown{E}\le 0$. Then there exsits a birational $G$-model $(S',G)$ of $(S,G)$ satisfying the following $2$ conditions
\begin{enumerate}
\item[$\bullet$] \noindent  $(S',G)$ is an action with {\rm IFP}, and
\item[$\bullet$] \noindent $\kappa (K_{S'})=-\infty$.
\end{enumerate} 
\end{lemma}

\begin{proof}
We can assume $\bar{S}=S$. Thanks to the argument in Section $2$, from $(S,G)$, we can construct an action $(S',G)$ with IFP, satisfying if $(S^*,G)$ is a common resolution of $S$ and $S'$,  
$$
\xymatrix{
 & S^*\ar[d]^{\pi }\ar[r]^{\phi } & S \ar@{-->}[dl]^f \\
     &                  S'                                &
}$$ 
the exceptional divisors of $\pi$ only consists of curves in $\sum_S$ and the exceptional curves of $\phi$.  

Let $\pi^*(K_{S'})=K_{S^*}+ E+ F$, where $\supp (E)$ is in the birational transform of $\sum_S$, $F$ is exceptional for $\phi$, then from $\rdown{E+F} \le 0$, we know $\rdown{E} \le 0$. Hence for any $m\in \N$,   
$$h^0(mK_{S'})=h^0(m(K_{S^*}+ E+ F)) \le h^0(m\phi_*(K_{S^*}+ E +F))=h^0(m(K_{S}+ \phi_*E)=0.$$

\end{proof}

\begin{lemma}
 The action $(\P^2,G)$ as in (\ref{example}) is not birational to any action on a log del Pezzo surface $(S,G)$ with IFP.
\end{lemma}
\begin{proof}
By way of contradiction, if there is such a surface $S'$, its minimal resolution $\pi:S^* \to S'$ has a equivariant morphism $\phi$ to $\P^2$, and $\pi$ contracts the birational transform of the above lines. 
$$
\xymatrix{
 & S^*\ar[d]^{\pi }\ar[r]^{\phi } & \P^2 \ar@{-->}[dl]^f \\
     &                  S'                                &
}$$

Since $S'$ contains only quotient singularities, in the exceptional locus of the morphism $\pi: S^*\to S'$, any 3 irreducible component cannot intersect at an identical point, which implies any 3 components of the birational transform of $\sum_{\P^2}$ on $S^*$ can not intersect at an identical point. Hence the morphism $\phi:S^*\to \P^2$ must factor through the surface which we achieve by blowing up the 12 points $\{ (1,\epsilon_3^i,\epsilon_3^j)(0\le i,j\le 2), (1,0,0),(0,1,0),(0,0,1)\} $ on $\P^2$. Then to get $S^*$ from $\P^2$, we have to blow up at least 4 points on each one of the 9 lines $\{x_i=\epsilon_3^k x_j\}$. Thus the self-intersection numbers of the birational transform of these lines on $S^*$ are smaller or equal to $-3$. If we write $K_{S^*}+\sum_i a_iE_i=K_S$ with $0<a_i<1$, then the coefficients of the above lines are greater or equal to $\frac{1}{3}$ (cf. \cite{alexeev}, 2.17). Hence for $m$ divisible enough
$$H^0(S,mK_S)=H^0(m(K_{S^*}+\sum_i a_iE_i))=H^0(\P^2,m\phi_*(K_{S^*}+\sum_i a_iE_i)),$$
The second equality holds because $K_{S^*}+\sum_i a_iE_i-\phi^*\phi_*(K_{S^*}+\sum_i a_iE_i) \ge 0$, which is implied by the nefness of $-(K_{S^*}+\sum_i a_iE_i)$ (cf. \cite{kollarmori}, 3.39). But thanks to the computation of $a_i$, we know  $\phi_*(K_{S^*}+\sum_i a_iE_i)\ge K_{\P^2}+1/3\sum_{i=1}^{9}L_i \ge 0$, here $L_i$'s means the divisors of the above 9 lines.

\end{proof}

\begin{proof}[Proof of (\ref{main2})] By (\ref{mmp}), it suffices to prove that given an action $(S,G)$ as in (1)-(4) of (\ref{main}), we can choose $(\bar{S},G)$ in the same $G$-birational class satisfying the conditions of (\ref{kodaira}). According to the above example, we know groups containing $(\Z/3)^2:\Z/2$ cannot act on any log Del Pezzo surfaces with IFP. We check the remaining cases in section 3 as follows:

For (\ref{intransitive}), after blowing up the origin, we get a ruled surface $\F_1$, and $\sum_{\F^1}$ consists of two sections and a set of fibers. Denote the fiber class as $L$. Then $(K_{\F_1}+E)\cdot L <0$, as the coefficient of the sections in $E$ are less than $1$. Since $|L|$ is a covering family, $K_{\F_1}+E$ is not pseudo-effective.

For (\ref{transimpri}), $\Z/n : \Z/3$ acts on $\P^2$ with IFP. For $S_3$, blow up the intersection point $(1,1,1)$ of curves in $\sum_{\P^2}$, then we can argue as in the above case.

For (\ref{F_0}), if $G$ is $G_1\times G_2$, then one of them, say $G_1$, is cyclic ($|G_1|, |G_2|$ are coprime). So $\sum_{\F^0,G}$ contains at most $2$ sections for the corresponding fiberation. Denote the class of the fiber as $L$, we have $(K_{\F_0}+E)\cdot L <0$, which impies $K_{\F_0}+E$ is not pseudo-effective. For general $G=(G_1,H_1,G_2,H_2)_{\alpha}$, we can argue in the same way because $\sum_{\F^0,G}=\sum_{\F^0,H_1\times H_2}$. If $G : G^0 = 2$, the only new groups are $F_{4n}$, $G_{4n}$ and $H_{4n}$. For these cases, $\sum_G$ are always empty. 

For (\ref{degree6}), by the proof there, we know for $\Z/2\times (\Z/n:\Z/3)$, the surface $S$ itself gives a model of action with IFP. For $S_3$, it is equivariantly birational to an action on $\P^2$ (cf. \cite{dol}, 8.1).

For (\ref{degree5}), as we discussed before, it is birational to an action on $\F_0$.
\end{proof}

\section{$\pi_1$ of smooth loci of log del pezzo surfaces}

In the previous section, we give a table containing precisely $\pi_1$ of smooth points of log del Pezzo surfaces and all their quotient groups. In this section, we aim to determine which of these groups can be actual fundalmental groups. In other words, for a given $G$, we want to construct a log del Pezzo surface $S$ such that $G$ acts on it with IFP and $\pi_1(S^{sm})={e}$.

\begin{prop}
Every group $G$ in (3), (4) of (\ref{main}) is $\pi_1$ of smooth points of some log Del Pezzo surface.
\end{prop}
\begin{proof}
We observe that every group $G$ in (3), (4) of (\ref{main}) has the property: for the action $(S,G)$ arises from the classification in Section 3, $\sum_{S,G}$ contains at most one irreducible curve and $\pi_1(S \backslash \sum_{S,G})=\{ e \}$. Blow up a general orbit of $G$ on $\sum_{S,G}$ and then contract its birational transform, we get an equivariant model $(\tilde{S},G)$ such that:
\begin{enumerate} 
\item $G$ acts on $\tilde{S}$ with IFP, 
\item $\tilde{S}$ is a log Del Pezzo surface,
\item $\tilde{S}^{sm}$ contains $S \backslash \sum_{S,G}$ as an open set.
\end{enumerate}
So we conclude $\pi_1(\tilde{S}^{sm})=\{e\}$ and $G$ can be $\pi_1$ of smooth points of some log Del Pezzo surface. 
\end{proof}
The remaining cases are subgroups of $\GL_2(\C)$ or $\PGL_2(\C) \times \PGL_2(\C)$.

\begin{prop}\label{birationalF_e}
Given a group $G$ in (\ref{intransitive}) (resp. (\ref{F_0})), if $\tilde{G}$ (resp. $G=G'$) has the form $(\mu_{mk},\mu_m,G_1, G_2)$, then it is a fundamental group of some log del Pezzo surface. 
\end{prop}
\begin{proof}
If $m=1$, then $G$ is either polyhedral or binary polyhedral.  For any binary polyhedral group $G$, consider its action on $\P^2$ which factors through $\SL_2(\C)$. The only possible component of $\sum_{\P^2,G}$ is the infinite line $L$.  Blow up a $G$-orbit on $L$, and then contract $L$, we get a pair $(S,G)$ with IFP. Since $S^{sm}$ contains $\C^2$ as an open set, $\pi_1(S^{sm})=\{e\}$.  

From now on we assume $m \not= 1$. In the case (\ref{intransitive}), we blow up the original point, and assume we always have $S=\F_e$. Then the configuration of $\sum_{S,F}$ is as following,

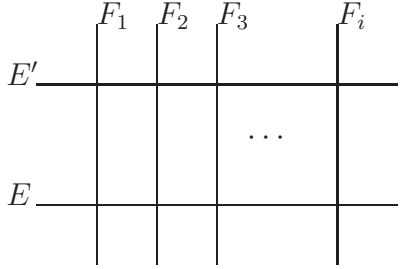
\begin{center}
\setlength{\unitlength}{1mm}
\begin{picture}(40,40)
\multiput(8,0)(8,0){3}{\line(0,1){32}}
\multiput(0,8)(0,16){2}{\line(1,0){48}}
\put(28,16){$\cdots$}
\put(40,0){\line(0,1){32}}
\put(8,32){$F_1$}
\put(16,32){$F_2$}
\put(24,32){$F_3$}
\put(40,32){$F_i$}
\put(-4,8){$E$}
\put(-4,24){$E'$}
\end{picture}
\end{center} 

\noindent where the vertical lines are fibers of $\F_e$. We will do the following sequence of $G$-birational operations on $\F_e$, which terminates with a $G$-surface $S$ satisfying $\pi_1(S^{sm})={e}$.\\ 
 
\noindent Step(1): First We construct a birational model $(S,G)$ such that $\sum_{S,G}$ does not contain any vertical lines. 

The way to construct $S$ is as follows: assume $E^2\le 0$, we first equivariantly blow up the intersection points of $F_i's$ and $E$, then contract $F_i's$. By (\ref{separation}), we know after finite steps of such operation at each intersection point, we will have a $G$-surface $S=\F_{r}$ such that $\sum_{S,G}$ does not contain any fiber. \\

\noindent Step(2): We construct a model $S$ as in step(1) with the additional property that $-q\le E^2<  0 $, where $q$ is the length of the $G$-orbit of an general point on $E'$

The way to construct $S$ is similar as in step(1). Assume $E^2=-r<0, (E')^2=r$. If we choose a general point $x$ on $E'$, so its stabilizer $G_x$ is precisely the subgroup whose elements fix $E'$ pointwise. We conclude $q=|G|/|G_x|$. Now blow up these $q$ points, and contract the birational transforms of the fibers which pass through them. We have a new ruled surface with $E^2=-r+q$, $(E')^2=r-q$. By the generality of the $q$ points on $E'$, we know $G$ acts on this new surface with $\sum=\{E,E'\}$. 

\noindent step(3): We have to deal with 2 different cases.

subcase(1): if $-q<E^2$, first we blow up a generic orbit on $E'$, then contract $E$ and $E'$ (the contractabilty of $E'$ comes from the assumption $-q<E^2$). The resulting surface $S$ is a log Del Pezzo surface which $G$ acts with IFP (see the last part of step(2)) and $\Pic(S)^G=\Z$. Then apply the computation in \cite{mumford}, we can easily conclude that $\pi_1(S^{sm})=\{e \}$.

subcase(2): if $E^2=-q$. We again start with choosing a general orbit of $q$ points on $E'$, and assume the fibers passing through the points are $F_j (1\le j \le q)$. We blow up these $q$ points, with exceptional divisors $E^{(1)}_j $; then blow up the intersection points of all $E^{(1)}_j$ with $E'$, with exceptional divisors $E^{(2)}_j$; then again blow up the intersection points of all $E^{(2)}_j$ with $E'$, etc. We do this type of blow-up $p>1$ times, where $p$ satisfies $\gcd(q,p)=1$.\\

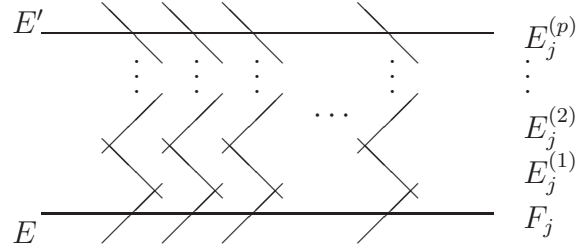
\begin{center}
\setlength{\unitlength}{1mm}
\begin{picture}(60,30)
\multiput(8,0)(8,0){3}{\line(1,1){8}}
\multiput(16,6)(8,0){3}{\line(-1,1){8}}
\multiput(8,12)(8,0){3}{\line(1,1){8}}
\multiput(16,24)(8,0){3}{\line(-1,1){8}}
\multiput(12,20)(8,0){3}{$.$}
\multiput(12,22)(8,0){3}{$.$}
\multiput(12,24)(8,0){3}{$.$}
\multiput(0,4)(0,24){2}{\line(1,0){60}}
\put(36,16){$\cdots$}
\put(42,0){\line(1,1){8}}
\put(50,6){\line(-1,1){8}}
\put(42,12){\line(1,1){8}}
\put(50,24){\line(-1,1){8}}
\put(46,24){$.$}
\put(46,20){$.$}
\put(46,22){$.$}
\put(-4,0){$E$}
\put(-4,28){$E'$}
\put(64,2){$F_j$}
\put(64,8){$E_j^{(1)}$}
\put(64,14){$E_j^{(2)}$}
\put(64,24){$.$}
\put(64,20){$.$}
\put(64,22){$.$}
\put(64,26){$E_j^{(p)}$}
\end{picture}
\end{center}

\noindent Now we contract the curves $E$, $E'$ and $E^{(1)}_j$, $E^{(2)}_j$,..., $E^{(p-1)}_j$, so we have the demanding log Del Pezzo surface $S$ satisfying $\Pic(S)^G=\Z$. In fact, the only element of $G$ which fixes $E^{(i)}_j$ pointwise is $e$, thus $G$ acts $S$ with IFP. By \cite{mumford}, we know $\pi_1(S^{sm})$ is the finite group generated by $a,b,c_j (1\le j\le q)$ with the relations $a^p=1, b^{pq-q}=1, (c_j)^p=1, a=c_j^{-1},  b=c_j$. But $\gcd(p,q)=1$ implies this group is in fact trivial. \\ 
\end{proof}

\begin{rem}
The remaining cases are when $S=\P^1 \times \P^1$, $G$ is one of the following groups:
$(D_{2m},\Z/m,O,T)$ ($\gcd(m,6)=1$), $(D_{6m},\Z/m,O, (\Z/2)^2)$ $(\gcd(m, 2) =
1$ or $(D_{2m},\Z/m,D_{4n},D_{2n})$ ($\gcd(m,2n)=1$). Arguing in a similar way as above, we know these groups can act on a del Pezzo surface $S$ with IFP and $\pi_1(S^{sm})=\Z/2$. In fact, this seems to be the best thing we can get out of these cases, namely,
\begin{q}
For any (rank 1) log del Pezzo surface $S$ with a finite group $G$ acting on it, such that the action $(S,G)$ is with {\rm IFP} and  birational to one of the abvoe actions, then $\Z/2\subset \pi_1(S^{sm}).$
\end{q}
\end{rem}

\end{document}